\documentstyle[12pt]{article}

\catcode`\@=11
\begin{document}

\centerline{\large\bf SUPER YANGIAN DOUBLE $DY\left( gl(1|1)\right)$}
\centerline{\large\bf AND ITS GAUSS DECOMPOSITION }
\vspace{0.8cm}
\centerline{\sf Jin-fang Cai, Guo-xing Ju, Ke Wu }
\baselineskip=13pt
\centerline{ Institute of Theoretical Physics, Academia Sinica, }
\baselineskip=12pt
\centerline{ Beijing, 100080, P. R. China }
\vspace{0.3cm}
\centerline{ and}
\vspace{0.3cm}
\centerline{\sf Shi-kung Wang}
\baselineskip=13pt
\centerline{  CCAST (World Laboratory), P.O. Box 3730, Beijing, 100080, 
  P. R. China  }
\baselineskip=12pt
\centerline{ and  Institute of Applied Mathematics, Academia Sinica,}
\baselineskip=12pt
\centerline{Beijing, 100080, P. R. China }
\vspace{0.9cm}
\begin{abstract} {We extend Yangian double to super (or graded) case 
and give its
  Drinfel$^{\prime}$d generators realization  by Gauss decomposition.}
\end{abstract}

\baselineskip=14pt
\vspace{1cm}

Quantum algebras(in general meaning, it includes Yangians and quantum
affine algebras etc.) are new algebraic structures discovered
 about ten years ago [1-4], they play important roles in the 
study of soluble statistical models and quantum fields theory.  
  Quantum universal enveloping algebras of  simple Lie algebras,
 which are related with some simple solutions (without spectral parameter)
of Yang-Baxter equation, have been extensively and
 deeply studied in the past few years.
Quantum affine algebras and Yangians are related respectively with
 trigonometric and rational solutions of Yang-Baxter equation,
they have three realizatons in literatures: Chevalley generators,
$T^{\pm}$-matrix and Drinfel$^{\prime}$d generators. The first realization
was proposed independently by Drinfel$^{\prime}$d and Jimbo, the second 
realization  has direct meaning in quantum inverse scattering method and 
it's convenient to introduce central extension using this realization 
[5]. The isomorphism of $T^{\pm}$-matrix and  Drinfel$^{\prime}$d generators
realizations of quantum affine algebra was established through 
Gauss decomposition by Ding and Frenkel [17] .
  The Yangian can be viewed as a deformation of only half
 of the corresponding loop algebra , while Yangian double are 
deformation of the complete loop algebra. 
 Yangian doubles (and with central extension) have been studied by
Bernard, Khoroshkin and Iohara etc.[6-12]. The properties of
 super Yangians and its 
representations have been studied by some authors [13,14]. 
But to the author's knowledge, 
the super Yangian doubles and Gauss decomposition in super case have not been 
studied up to now . In this paper, we extend Yangian double to  super
 case and give its 
Drinfel$^{\prime}$d generators realization through Gauss decomposition. 

\vskip 1cm

Given a 2-dimensional graded vector space $V$ and let its second basis 
 be odd, then  the super (graded) Yang-Baxter equation [15-17]
takes this form:
 $$ 
\eta _{12}R_{12}(u) \eta _{13}R_{13}(u+v) \eta _{23}R_{23}(v) 
= \eta _{23}R_{23}(v) \eta _{13}R_{13}(u+v) \eta _{12}R_{12}(u),
\eqno (1)
 $$
where $R(u) \in End(V\otimes V) $ and it must obey weight conservation 
condition: $R_{ij,kl}\neq 0 $ only when $i+j = k+l $.
 $\eta _{ik,jl}=(-1)^{(i-1)(k-1)} \delta _{ij} \delta _{lk}$ .
In dealing with tensor product in graded case,
 we must use this form:
$(A\otimes B)(C \otimes D)=(-1)^{P(B)P(C)}AC \otimes BD $.
The super Yang-Baxter equation can also be writed in components as follow:
$$
\begin{array}{rcl}
&&R_{ib,pr}(u)R_{pa,js}(u+v)R_{rs,dc}(v)(-1)^{(r-1)(s+a)} \\
&=&
(-1)^{(e-1)(f+c)}R_{ba,ef}(v)R_{if,kc}(u+v)R_{ke,jd}(u)
\end{array} \eqno (2)
$$

It is  easily to prove that $R(u)=u I+\hbar{\cal P} $ satisfy the above 
super Yang-Baxter equation, here ${\cal P}_{ik,jl}=(P\eta)_{ik,jl}=
(-1)^{(i-1)(k-1)} \delta _{il} \delta _{jk} $ is the 
permutation operator in graded case. 
   
In similiar to $DY(gl_2) $ [10],  Super Yangian Double $DY\left( gl(1|1)
 \right)$ is a Hopf algebra generated 
by $ \{ t_{ij}^k \vert 1\leq i,j \leq 2, k\in {\bf Z }\} $ which obey the 
following relations:
$$
\begin{array}{rcl}
 R(u-v)T_1 ^{\pm}(u) \eta T_2 ^{\pm}(v) \eta &= &\eta T_2 ^{\pm}(v) \eta T_1 ^{\pm}(u) R(u-v) \\
R(u-v)T_1 ^{+}(u) \eta T_2 ^{-}(v) \eta &= &\eta T_2 ^{-}(v) \eta T_1 ^{+}(u) R(u-v) 
\end{array}  \eqno (3)
$$
here we use standard notation:
$T_1^{\pm}(u)=T^{\pm}(u)\otimes {\bf 1}, T_2^{\pm}(u)={\bf 1} \otimes
T^{\pm}(u) $ and $T^{\pm}(u)=\left( t_{ij}^{\pm}(u) \right)
 _{1 \leq i,j \leq 2}$,
  $t_{ij}^{\pm}(u)$ are generate functions of $t_{ij}^k$ :
$$ t_{ij}^{+}(u)=\delta _{ij}-\hbar\sum_{k\geq 0}t_{ij}^k u^{-k-1},
  {\hspace{0.5cm}}
 t_{ij}^{-}(u)=\delta _{ij}+\hbar\sum_{k < 0}t_{ij}^k u^{-k-1}.
\eqno (4)
$$

The Hopf pairing relation between $T^+(u)$ and $ T^-(v)$ is defined as:
$$< T^+_1(u) ~,~ T^-_2(v) > = R(u-v)  \eqno (5)$$.

From relations (3), we can get the (anti-)commutation relations
among $\{ t_{ij}^k(u) \}$:
$$
(u-v)[t_{ij}^{\sigma}(u) ~,~ t_{kl}^{\rho}(v)\}+(-1)^{ij+jk+ki+1}
\hbar (t_{kj}^{\sigma}(u)t_{il}^{\rho}(v)-t_{kj}^{\rho}(v)
t_{il}^{\sigma}(u))=0  \eqno (6)
$$
 where ($\sigma$ , $\rho$) take (+,+), ($-,-$) and ($+,-$). 
The parity of $t_{ij}^{\pm}$ equal to $(i-1)(j-1)$ which means that
($t_{11}^{\pm}(u), t_{22}^{\pm}(u) $ ) are even
 and ($t_{12}^{\pm}$ , $t_{21}^{\pm}(u)$) are odd.
We denote  $[  ~,~\}$ as super commutator which is anti-commutator
only when each of the elements in it is odd(or Fermionic) type. The 
follows are some special examples of the above relation:
$$
\begin{array}{ll}
&(u-v)[t_{11}^{\sigma}(u) ~,~ t_{12}^{\rho}(v)]+
\hbar(t_{11}^{\sigma}(u)t_{12}^{\rho}(v)-t_{11}^{\rho}(v)
t_{12}^{\sigma}(u))=0  \\
&(u-v)[t_{22}^{\sigma}(u) ~,~ t_{12}^{\rho}(v)]-
\hbar(t_{12}^{\sigma}(u)t_{22}^{\rho}(v)-t_{12}^{\rho}(v)
t_{22}^{\sigma}(u))=0  \\
&(u-v)\{t_{12}^{\sigma}(u) ~,~ t_{21}^{\rho}(v)\}-
\hbar(t_{22}^{\sigma}(u)t_{11}^{\rho}(v)-t_{22}^{\rho}(v)
t_{11}^{\sigma}(u))=0  
\end{array} \eqno (7)
$$

 Hopf structure of $DY\left(gl(1\vert 1)\right)$  are defined
as follows:

$$
\begin{array}{ll}
&\triangle t_{ij}^{\pm}(u) =\sum_{k=1,2} t_{kj}^{\pm}(u) \otimes 
t_{ik}^{\pm}(u)(-1)^{(k+i)(k+j)}, \\
& \epsilon\left(t_{ij}^{\pm}(u)\right)=\delta _{ij} ,{\hspace{1.cm}}
S\left( ^{st} T^{\pm}(u) \right) =\left[ ^{st} T^{\pm}(u) \right]^{-1}.
\end{array} \eqno (8)
$$
here $\left[^{st} T^{\pm}(u) \right]_{ij} =(-1)^{i+j}t_{ji}^{\pm}(u) $ .

We can get the Drinfel$^{\prime}$d generators realization of $DY(gl(1|1))$ 
by Gauss decomposition in the same way which has been used in $DY(gl_2)
$. Introducing  transformation for generate functions as follows:
$$  
E^{\pm}(u) =\frac{1}{\hbar}t_{11}^{\pm}(u+\frac{\hbar}{2})^{-1}
t_{12}^{\pm}(u+\frac{\hbar}{2}),{\hspace{0.5cm}}
F^{\pm}(u)=\frac{1}{\hbar}t_{21}^{\pm}(u+\frac{\hbar}{2})t_{11}^{\pm}
(u+\frac{\hbar}{2})^{-1} , 
$$
$$
k_1^{\pm}(u) =t_{11}^{\pm}(u), {\hspace{0.5cm}} k_2^{\pm}(u) = t_{22}^{\pm}(u)-
t_{21}^{\pm}(u)t_{11}^{\pm}(u)^{-1}
t_{12}^{\pm}(u) . \eqno (9)
$$
it means 
$$
\begin{array}{rcl}
T^{\pm}&=&\left(\begin{array}{cc}
    1              &\\
 \hbar F^{\pm}(u-\frac{\hbar}{2}) &1
\end{array}\right)
\left(\begin{array}{cc}
   k_1^{\pm}(u)             &\\
  & k_2^{\pm}(u)
\end{array}\right)
\left(\begin{array}{cc}
    1   & \hbar E^{\pm}(u-\frac{\hbar}{2})\\
   & 1
\end{array}\right)  \\
&=& \left( \begin{array}{cc}
k_1^{\pm}(u) & \hbar k_{1}^{\pm}(u)E^{\pm}(u-\frac{\hbar}{2})\\
\hbar F^{\pm}(u-\frac{\hbar}{2}) k_1^{\pm}(u) & k_2^{\pm}(u)+\hbar^2
F^{\pm}(u-\frac{\hbar}{2}) k_1^{\pm}(u)E^{\pm}(u-\frac{\hbar}{2})
\end{array} \right) 
\end{array} \eqno (10)
$$
and let
$$ 
E(u)=E^+(u)-E^-(u), {\hspace{0.5cm}} F(u)=F^+(u)-F^-(u) ,
$$
$$
H^{\pm}(u)=k_2^{\pm}(u+\frac{\hbar}{2})k_1^{\pm}(u+\frac{\hbar}{2})^{-1},
 {\hspace{0.5cm}}
K^{\pm}(u)=k_2^{\pm}(u+\frac{\hbar}{2})k_1^{\pm}(u-\frac{\hbar}{2}) .
\eqno (11)
$$
By the similiar calculate process which have been used in quantum 
affine algebra [18],
 we can obtain that these new generate 
functions $E(u),F(u), H^{\pm}(u) $ and $K^{\pm}(u)$ 
satisfy the following commutation relations:
$$
\begin{array}{ll}
&[H^{\sigma}(u)~,~H^{\rho}(v)]=[H^{\sigma}(u)~,~K^{\rho}(v)]~,~ 
\forall \sigma,\rho=+,- \\
&[K^{\sigma}(u)~,~K^{\rho}(v)]=0,  \\
&[H^{\pm}(u)~,~E(v)]=[H^{\pm}(u)~,~ F(v)]=0, \\
&E(v)K^{\pm}(u)=\frac{u-v+\hbar}{u-v-\hbar}K^{\pm}(u)E(v), \\
&F(v)K^{\pm}(u)=\frac{u-v-\hbar}{u-v+\hbar}K^{\pm}(u)F(v) , \\
&\{ E(u)~,~E(v)\}=\{F(u)~,~F(v)\}=0 ,\\
&\{E(u)~,~F(v)\}=\frac{1}{\hbar}\delta(u-v)(H^-(v)-H^+(u)). 
\end{array} \eqno (12)
$$

The  coproduct and counit structure for $E(u), F(u) ,H^{\pm}(u) $ and $K^{\pm}(u)$
can be calculate out as follows:
$$
\begin{array}{ll}
&\triangle E^{\pm}(u)= E^{\pm}(u)\otimes 1+ K^{\pm}(u) \otimes E^{\pm}(u) \\
&\triangle F^{\pm}(u)=1\otimes F^{\pm}(u)+F^{\pm}(u)\otimes k^{\pm}(u) \\
&\triangle K^{\pm}(u)= K^{\pm}(u)\otimes K^{\pm}(u)-
 F^{\pm}(u-\hbar) K^{\pm}(u)\otimes (E^{\pm}(u)K(u)+K(u)E^{\pm}(u)) \\
&\triangle H^{\pm}(u)=H^{\pm}(u)\otimes H^{\pm}(u) \\
&\epsilon (K^{\pm}(u)) =\epsilon (H^{\pm}(u)) =1 \\
&\epsilon (E^{\pm}(u)) =\epsilon (F^{\pm}(u)) =0  
\end{array}   \eqno (13)
$$

Expanding $H^{\pm}(u),K^{\pm}(u),E(u),F(u) $ in Fourier series
$$
\begin{array}{ll}
&E(u)=\sum_{k\in Z}e_{k}u^{-k-1},~
  F(u)=\sum_{k\in Z}f_{k}u^{-k-1}, \\
& H^{+}(u)=1+\hbar\sum_{k\ge 0}h_{k}u^{-k-1},~
  H^{-}(u)=1-\hbar\sum_{k<0}h_{k}u^{-k-1}, \\
& K^{+}(u)=1+\hbar\sum_{l\ge 0}k_{l}u^{-l-1},~
  K^{-}(u)=1-\hbar\sum_{l<0}k_{l}u^{-l-1}.
\end{array} \eqno (14)
$$
  then the commutation relations among $h_l, k_l, e_l $ and $f_l$
 give  Drinfeld's generator realization of $DY(gl(1|1))$ :

$$
\begin{array}{ll}
& [h_{k},h_{l}]=[k_k, k_l]=[h_k,k_l]=0 \\
&[h_k,e_l]=[h_k,f_l] =0 \\
&[k_0,e_l]=-2e_l~,~[k_0,f_l]=2f_l \\
&[k_{k+1},e_l]-[k_k,e_{l+1}]+\hbar(k_k e_l+e_l k_k)=0 \\
&[k_{k+1},f_l]-[k_k,f_{l+1}]-\hbar(k_k f_l+f_l k_k)=0 \\
&\{e_k, f_l\}=-2h_{k+l}\\
&\{e_k, e_l\}=\{f_k, f_l\}=0 
\end{array}  \eqno (15)
$$

It is naturally that  maybe we can also discuss central  extension of 
 $DY(gl(1|1) $ and the  Gauss decomposition can also be applied to
 quantum affine superalgebra to  get Ding-Frenkel map of its two 
realizations.   As to its  physical applicant,
Yangian symmetry have been discovered in integrable quantum 
field [6,7], spin chains with long range interactions (including
Calogero-Sutherland models) [19,20], Hubbard model [21],
conformal field [22,23] and so on. Super Yangian has also been applied in 
color Calogero-Sutherland models [24] . Maybe the 
symmetries in these systems are the corresponding (super) Yangian doubles
 as ref.[6].

\vskip 1cm
{\bf \large References:}
\begin{enumerate}
 
\item[{[1]}] V.G. Drinfel$^{\prime}$d:
 ``Quantum Groups'', {\sl Proc.ICM-86(Berkeley)}. Vol.1. 
 NewYork Academic Press 1986, 789

\item[{[2]}]  V.G. Drinfel$^{\prime}$d, 
{\sl Soviet Math. Dokl.} {\bf 32} (1985) 254

\item[{[3]}] V.G. Drinfel$^{\prime}$d:
{\sl Soviet Math. Dokl.} {\bf 36} (1988) 212

\item[{[4]}] L.D. Faddeev, N.Yu. Reshetikhin and L.A. Takhtajan :
 {\sl Algebraic analysis } {\bf 1} (1988) 129

\item[{[5]}] N. Yu. Reshetikhin and M.A. Semenov-Tian-Shansky:
{\sl Lett. Math. Phys.} {\bf 19} (1990) 133

\item[{[6]}] D. Bernard and A. LeClair:  
      {\sl Nucl. Phys.} {\bf B399} (1993) 709

\item[{[7]}] A. LeClair and F. Smirnov: 
 {\sl Int. J. Mod. Phys. } {\bf A7} (1992) 2997

\item[{[8]}] S.M. Khoroshkin and V.N. Tolstoy: 
  {\sl Lett. Math. Phys.} {\bf 36} (1996) 373

\item[{[9]}] S.M. Khoroshkin:
 ``Central Extension of the Yangian Double '' , 
{\sl q-alg/9602031 }

\item[{[10]}]  S. Khoroshkin, D. Lebedev and S. Pakuliak:
`` Intertwining Operators for the Central
 Extension of the Yangian Double '',
 {\sl q-alg/9602030}

\item[{[11]}] K. Iohara and M. Kohno:
 {\sl Lett. Math. Phys. } {\bf 37} (1996) 319

\item[{[12]}] K. Iohara:
  {\sl J. Phys. } {\bf A29} (1996) 4593

\item[{[13]}] M.L. Nazarov:
  {\sl Lett. Math. Phys.} {\bf 21} (1991) 123

\item[{[14]}] R.B. Zhang :
 {\sl J. Math. Phys.} {\bf 36} (1995) 3854

\item[{[15]}] L. Liao and X.C. Song:
 {\sl Mod. Phys. Lett. } {\bf A6} (1991) 959

\item[{[16]}] T. Deguchi and Y. Akutsu :
 {\sl  J. Phys.} {\bf A23} (1990) 1861

\item[{[17]}] S. Majid and M.J. Rodr{\'i}gnez-Plaza:
 {\sl J. Math. Phys. } {\bf 36} (1995) 7081

\item[{[18]}] J. Ding and I.B. Frenkel: 
  {\sl Commun. Math. Phys.} {\bf 156} (1993) 277

\item[{[19]}] D. Bernard, M. Gaudin, F.D.M. Haldane and V. Pasquier:
 {\sl J. Phys.} {\bf A26} (1993) 5219

\item[{[20]}] F.D.M. Haldane, Z.N.C. Ha, J.C. Talstra, D. Bernard 
  and  V. Pasquier: 
    {\sl Phys. Rev. Lett.} {\bf 69 } (1992) 2021

\item[{[21]}]   D.B. Uglov and V.E. Korepin:
 `` The Yangian symmetry of the Hubbard Model''
{\sl hep-th/9310158}

\item[{[22]}]  K. Schoutens: 
   {\sl Phys. Lett.} {\bf B331} (1994) 335

\item[{[23]}] P. Bouwknegt, A.W.W. Ludwig and K. Schoutens:
  {\sl Phys. Lett.} {\bf B338} (1994) 448

\item[{[24]}] C. Ahn and W.M. Koo :
  {\sl  Phys. Lett.} {\bf B365} (1996) 105

\end{enumerate}

\end{document}